\documentclass{amsart}
\usepackage{amsmath,amsthm,amssymb,IMjournal}

\usepackage{amsmath}
\usepackage{amssymb}
\usepackage{graphicx}
\usepackage{color}

\newcommand{\noi}{\noindent}

\newcommand{\onh}{\frac{1}{2}}

\newcommand{\bc}[2]{\left.{#1}\right|_{#2}}

\newcommand{\nn}{\nonumber}
\newcommand{\bgeq}{\begin{equation}}
\newcommand{\eneq}{\end{equation}}
\newcommand{\bgeqr}{\begin{eqnarray}}
\newcommand{\eneqr}{\end{eqnarray}}

\newcommand{\leftcoll}{\begin{flushleft}\smallskip\nopagebreak\begin{minipage}[c]{12cm}\sloppy}
\newcommand{\midcoll}{\end{minipage}\kern 2cm\begin{minipage}[c]{12cm}\sloppy}
\newcommand{\rightcoll}{\end{minipage}\kern 0cm \smallskip\end{flushleft}}
\newcommand{\leftcolu}{\begin{flushleft}\smallskip\nopagebreak\begin{minipage}[c]{17cm}\sloppy}
\newcommand{\midcolu}{\end{minipage}\kern 2cm\begin{minipage}[c]{8cm}\sloppy}
\newcommand{\rightcolu}{\end{minipage}\kern 0cm \smallskip\end{flushleft}}
\newcommand{\leftcolo}{\begin{flushleft}\smallskip\nopagebreak\begin{minipage}[c]{18cm}\sloppy}
\newcommand{\midcolo}{\end{minipage}\kern 0.5cm\begin{minipage}[c]{9cm}\sloppy}
\newcommand{\rightcolo}{\end{minipage}\kern 0cm \smallskip\end{flushleft}}
\newcommand{\rc}{\textcolor{black}}

\newcommand{\bitleft}{\rm\leftmargini 1ex \begin{itemize}}
\newcommand{\bit}{\begin{itemize}}
\newcommand{\eit}{\end{itemize}}

\def\half{\frac{1}{2}}
\def\iXseg{\int_{u(X_{i-\half})}^{u(X_{i+\half})}\limits}

\begin{document}
\newtheorem{The}{Theorem}[section]

\numberwithin{equation}{section}

\title{Computational studies of conserved mean-curvature flow}

\author{\|Miroslav |Kol{\' a}{\v r}|, Prague,
        \|Michal |Bene{\v s}|, Prague,\\
        \|Daniel |{\v S}ev{\v c}ovi{\v c}|, Bratislava}

\rec {September 30, 2013}

\abstract The paper presents the results of numerical solution of the evolution law for the constrained mean-curvature flow. This law originates in the theory of phase transitions for crystalline materials and describes the evolution of closed embedded curves with constant enclosed area. It is reformulated by means of the direct method into the system of degenerate parabolic partial differential equations for the curve parametrization. This system is solved numerically and several computational studies are presented as well.
\endabstract




\section{Introduction}\label{sec1}

\noi
The article deals with the non-local mean-curvature flow described by the evolution law
\begin{eqnarray}\label{CMCFl}
v_{\Gamma} &=& - \kappa_{\Gamma} + \frac{1}{|\Gamma|}\int_\Gamma \kappa_{\Gamma} {\rm d}\rc{s},\\
\bc{\Gamma}{t=0} &=& \Gamma_{ini}, \nn
\end{eqnarray}
where $\Gamma$ is the closed curve in ${\mathbb R}^2$, ${\bf n}_\Gamma$ the normal vector to $\Gamma$, $v_\Gamma$ the velocity in the direction of the normal vector, $\kappa_\Gamma$ the (mean) curvature of $\Gamma$ and $F$ the external prescribed force.
\rc{Here $|\Gamma|$ is the length of $\Gamma$.}

Problem (\ref{CMCFl}) represents a variant of the mean curvature flow described as
\begin{eqnarray}\label{MCFl}
v_{\Gamma} &=& - \kappa_{\Gamma} + F,\\
\bc{\Gamma}{t=0} &=& \Gamma_{ini}, \nn
\end{eqnarray}
with a particular choice of the forcing term $F$, which is widely studied in the literature (see e.g. \cite{DziukDeck05}) as well as its various mathematical treatment by the direct (parametric) method (see e.g. \cite{BenJAP10,BenEJP09}), by the level-set method (see e.g. \cite{OsherSethian}) or by the phase-field method (see e.g. \cite{BenMCFl03}).

The constrained motion by mean curvature has been discussed in the literature as well (see \cite{Gage86,Dolce02,McCoy03,Ruuth08}).
In particular, problem (\ref{CMCFl}) has been mentioned in \cite{RuSte92,HenryHilMim11,MBNonlAC11} within the context of a modification of the Allen-Cahn equation
\cite{AllenCahn,CahnHil59(III)} approximating the mean-curvature flow \cite{BenMCFl03}.
The non-local character of the
equation is connected to the recrystallization phenomena where a fixed previously melted volume of the liquid phase solidifies again. It also can \rc{be applied} in dislocation dynamics in crystalline materials or in the digital image processing (see e.g. \cite{Sinica08}).
In this text, we treat (\ref{CMCFl}) by means of the parametric method and solve the resulting degenerate parabolic system numerically to provide the information on the solution behavior.

\begin{figure}\label{FVM}
\begin{center}
\includegraphics[width=\textwidth]{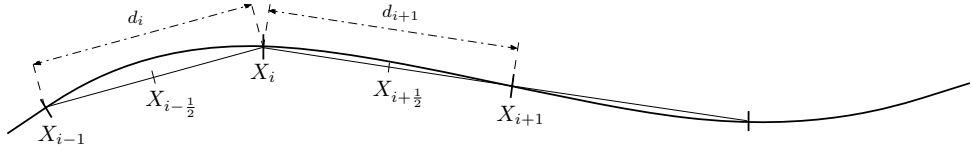}
\end{center}
\caption{Curve discretization by finite volumes.}
\end{figure}


\section{Equations}\label{sec2}

The direct method treating (\ref{CMCFl}) considers parametrization of the smooth time-dependent curve $\Gamma(t)$ by means of the mapping
$$\vec{X} = \vec{X}(t,u), \quad u \in S^1,$$
where $u$ is the parameter in a fixed interval. \rc{Here and after, we identify $S^1$ with the interval $[0,1]$ and we impose periodic boundary conditions on $\vec{X}$ at $u=0,1$.} Consequently the geometrical quantities of interest can be expressed by means of $\vec{X}$. The tangent vector and the normal vector are as follows
$${\bf t}_{\Gamma} = \frac{\partial_u \vec{X}}{|\partial_u \vec{X}|}, \quad {\bf n}_{\Gamma} = \frac{\partial_u \vec{X}^\perp}{|\partial_u \vec{X}|}.$$
The (mean) curvature is
\begin{equation}\label{Curv}
\kappa_{\Gamma} = -\frac{1}{|\partial_u \vec{X}|}\partial_u \left(\frac{\partial_u \vec{X}}{|\partial_u \vec{X}|}\right)\cdot{\bf n}_{\Gamma},
\end{equation}
and the normal velocity in the direction of ${\bf n}_{\Gamma}$ (the projection of the point velocity $\vec{v}_{\Gamma}$ at $\Gamma$ to ${\bf n}_{\Gamma}$) becomes
$$v_{\Gamma} = \vec{v}_{\Gamma} \cdot {\bf n}_{\Gamma}  \quad  \mbox{ where } \quad \vec{v}_{\Gamma} = \partial_t \vec{X}.$$
Substituting into (\ref{MCFl}) and assuming validity in the vectorial form yields the system
\begin{equation}\label{Dir}
\partial_t \vec{X} = \frac{1}{|\partial_u \vec{X}|}\partial_u \left(\frac{\partial_u \vec{X}}{|\partial_u \vec{X}|}\right)
+ F \frac{\partial_u \vec{X}^\perp}{|\partial_u \vec{X}|} \ \mbox{ in } \ (0,T)\times S^1
\end{equation}
known as the parametric (direct) description of (\ref{MCFl}).

Among advantages of this approach, an easy and straightforward treatment of the curve dynamics without additional approximation is offered. On the other hand, topological changes are not captured by it.

Further modifications of (\ref{Dir}) lead to the \rc{governing equation proposed by Dziuk et al.} in \cite{DziukDeck05} \rc{(see e.g. Bene{\v s} et al. \cite{BenEJP09} for applications in the dislocation dynamics)}
\begin{equation}\label{DirDeck}
\partial_t \vec{X} = \frac{\partial_{uu} \vec{X}}{|\partial_u \vec{X}|^2}
+ F \frac{\partial_u \vec{X}^\perp}{|\partial_u \vec{X}|} \ \mbox{ in } \ (0,T)\times S^1 .
\end{equation}
\rc{where}
\begin{equation}\label{DirNonl}
\rc{F} =  \frac{1}{\int_{S^1} |\partial_u \vec{X}| du}
\int_{S^1} \kappa_\Gamma(\vec{X}) |\partial_u \vec{X}| du
\end{equation}
with $\kappa_\Gamma(\vec{X})$ given by (\ref{Curv}), and the initial parametrization set as $\vec{X}|_{t=0} = \vec{X}_{ini}$.

\begin{figure}
\begin{center}
\includegraphics[width=0.48\textwidth]{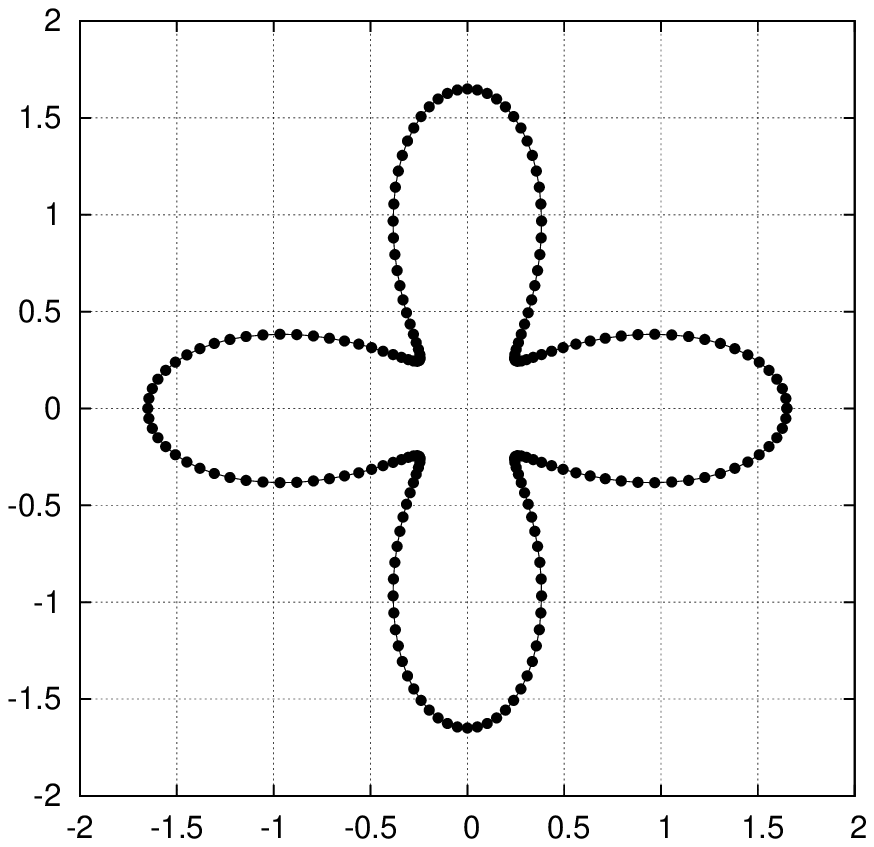}
\includegraphics[width=0.48\textwidth]{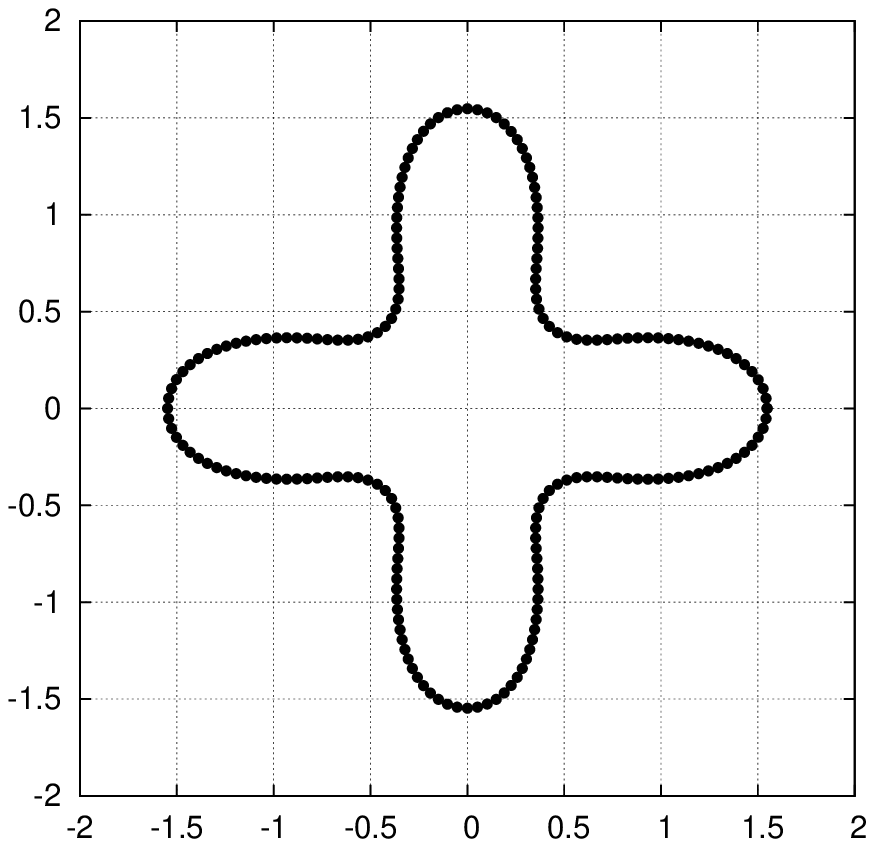}\\
\includegraphics[width=0.48\textwidth]{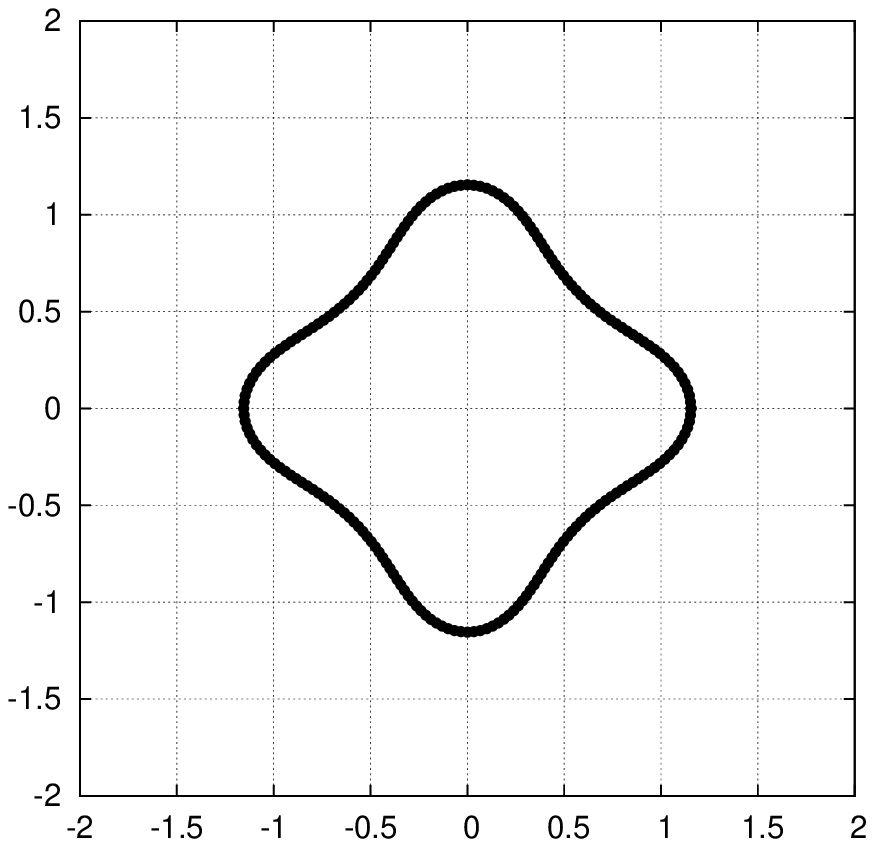}
\includegraphics[width=0.48\textwidth]{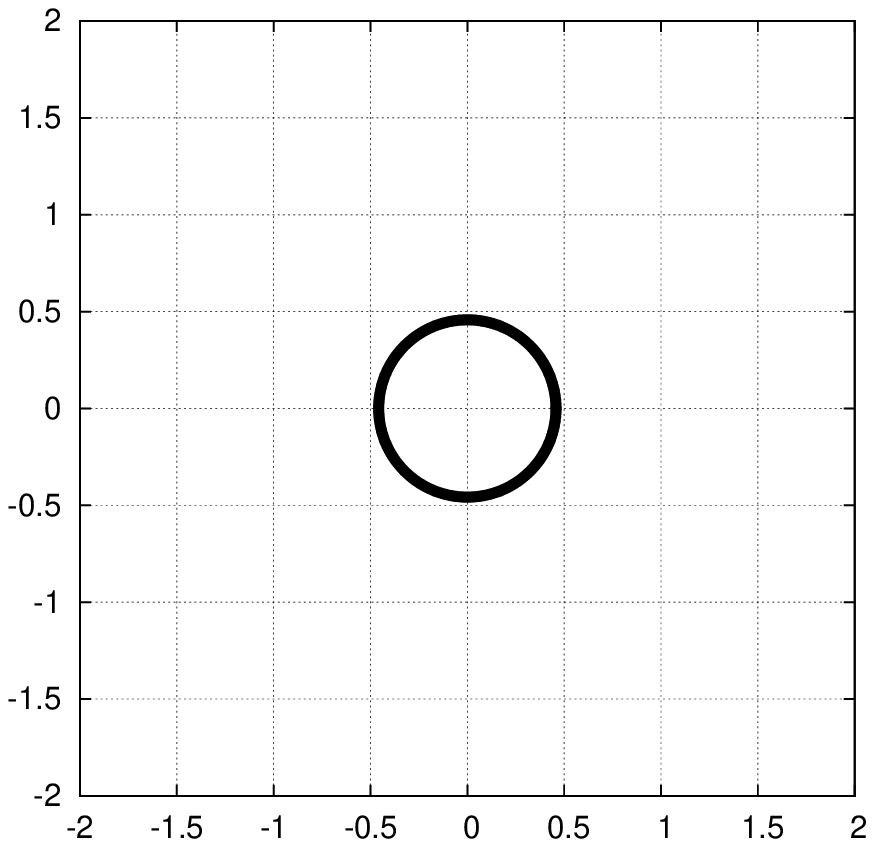}\\
\end{center}
\caption{Curve shortening flow (\ref{MCFl}) for which any closed curve shrinks to a point - the curve $\Gamma(t)$ is depicted for $t=0$, $t=0.025$, $t=0.125$ and $t=0.5$.}
\label{fig:ex1}
\end{figure}


\section{Numerical solution}\label{sec3}

For the discretization of (\ref{Dir}), the method of flowing finite volumes is used as e.g. in \cite{BenEJP09}. The discrete nodes $\vec{X}_i, i=0,\dots,M$\rc{,} are placed along $\Gamma(t)$ as shown in Figure \ref{FVM}.
The governing equation is integrated along the dual segments surrounding the nodes $\vec{X}_i, i=1,\dots,M-1$\rc{,}
\begin{eqnarray*}
  \iXseg \rc{\partial_t \vec{X}} |\partial_u \vec{X}| du &=& \iXseg \ \partial_u \left(\frac{\partial_u \vec{X}}{|\partial_u \vec{X}|}\right) du
    + \rc{F} \iXseg \partial_u \vec{X}^\perp du, \\
\rc{F} &=& \frac{1}{\int_{S^1} |\partial_u \vec{X}| du}
\int_{S^1} \kappa_\Gamma(\vec{X})|\partial_u \vec{X}| du , \\
\kappa_{\Gamma}(\vec{X}) &=& -\frac{1}{|\partial_u \vec{X}|}\partial_u \left(\frac{\partial_u \vec{X}}{|\partial_u \vec{X}|}\right)\cdot \frac{\partial_u \vec{X}^\perp}{|\partial_u \vec{X}|} .
\end{eqnarray*}
Resulting system of ordinary differential equations has the form
\begin{eqnarray}\label{eq:DCiScheme}
\frac{d\vec{X}_i}{dt}&=& \frac{2}{d_i + d_{i+1}}
\bigg(\frac{\vec{X}_{i+1} - \vec{X}_i}{d_{i+1}} - \frac{\vec{X}_i - \vec{X}_{i-1}}{d_i}\bigg)
+ F \frac{(\vec{X}^\bot_{i+1} - \vec{X}^\bot_{i-1})}{d_i + d_{i+1}} , \\
\kappa_i &=&  \frac{2}{d_i + d_{i+1}}
\bigg(\frac{\vec{X}_{i+1} - \vec{X}_i}{d_{i+1}} - \frac{\vec{X}_i - \vec{X}_{i-1}}{d_i}\bigg) \frac{(\vec{X}^\bot_{i+1} - \vec{X}^\bot_{i-1})}{d_i + d_{i+1}} , \nn \\
F &=& \frac{1}{\sum_{j=1}^M d_j} \sum_{j=1}^{\rc{M}} \kappa_j \frac{d_{j+1}+d_j}{2}, \nn \\
d_i &=& | \vec{X}_{i} - \vec{X}_{i-1} | \rc{, \quad d_{M+1} := d_1, \quad \vec{X}_0 := \vec{X}_M, \quad \vec{X}_{M+1} := \vec{X}_1} .  \nn
\end{eqnarray}
This system is solved by means of an semi-implicit backward Euler scheme. Details are similar to \cite{Sinica08}.

\begin{figure}
\begin{center}
\includegraphics[width=0.48\textwidth]{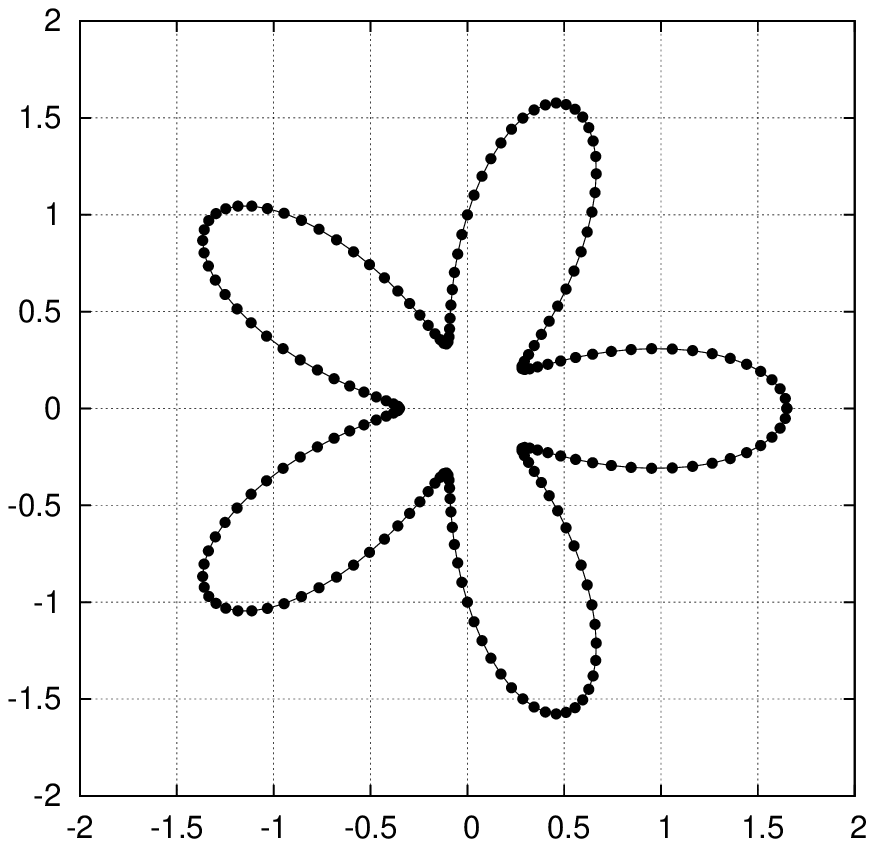}
\includegraphics[width=0.48\textwidth]{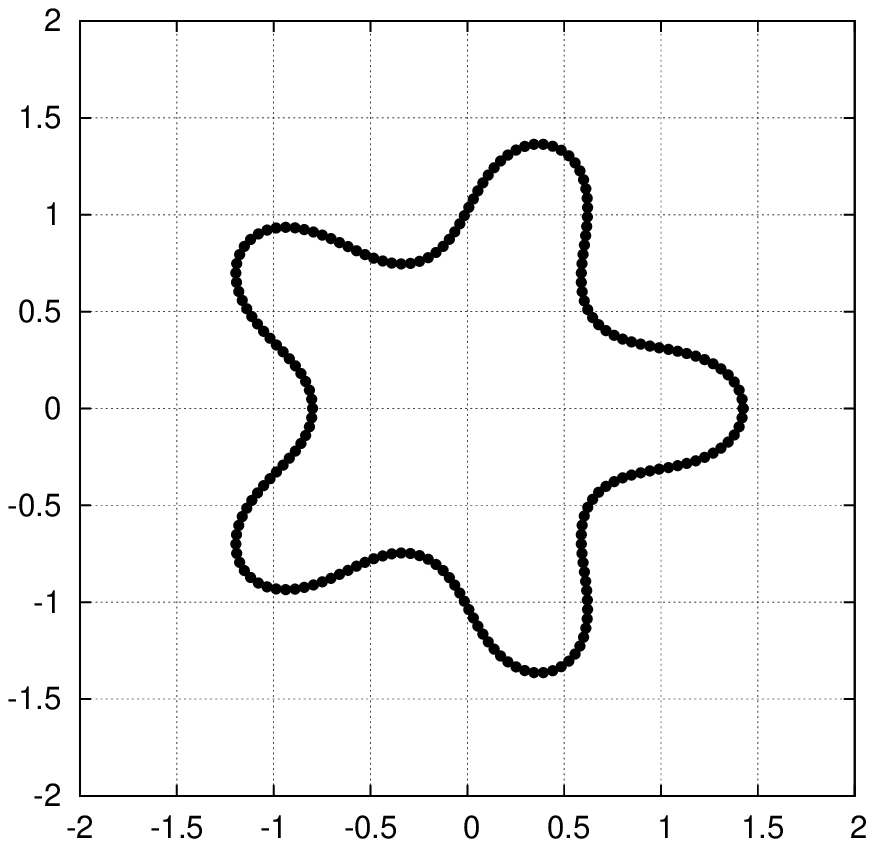}\\
\includegraphics[width=0.48\textwidth]{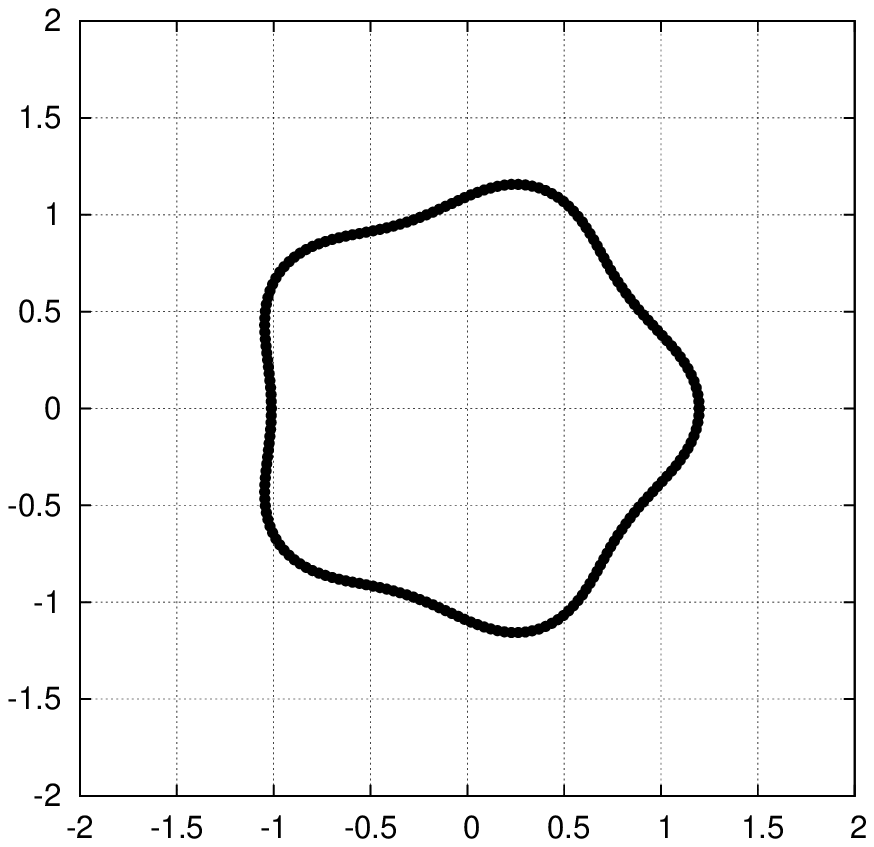}
\includegraphics[width=0.48\textwidth]{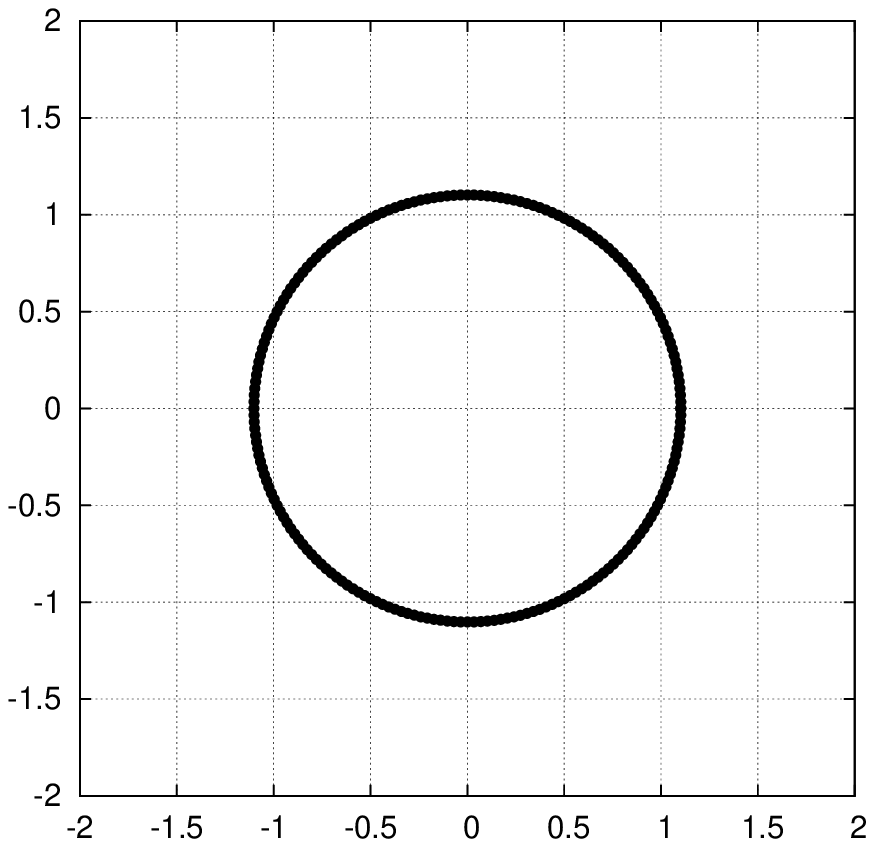}\\
\end{center}
\caption{\rc{Area-preserving} mean curvature flow (\ref{CMCFl}) where the initial 5-folded curve asymptotically attains the circular shape. The curve $\Gamma(t)$ is depicted for $t=0$, $t=0.05$, $t=0.125$ and $t=0.5$.}
\label{fig:ex2}
\end{figure}

\begin{figure}
\begin{center}
\includegraphics[width=0.48\textwidth]{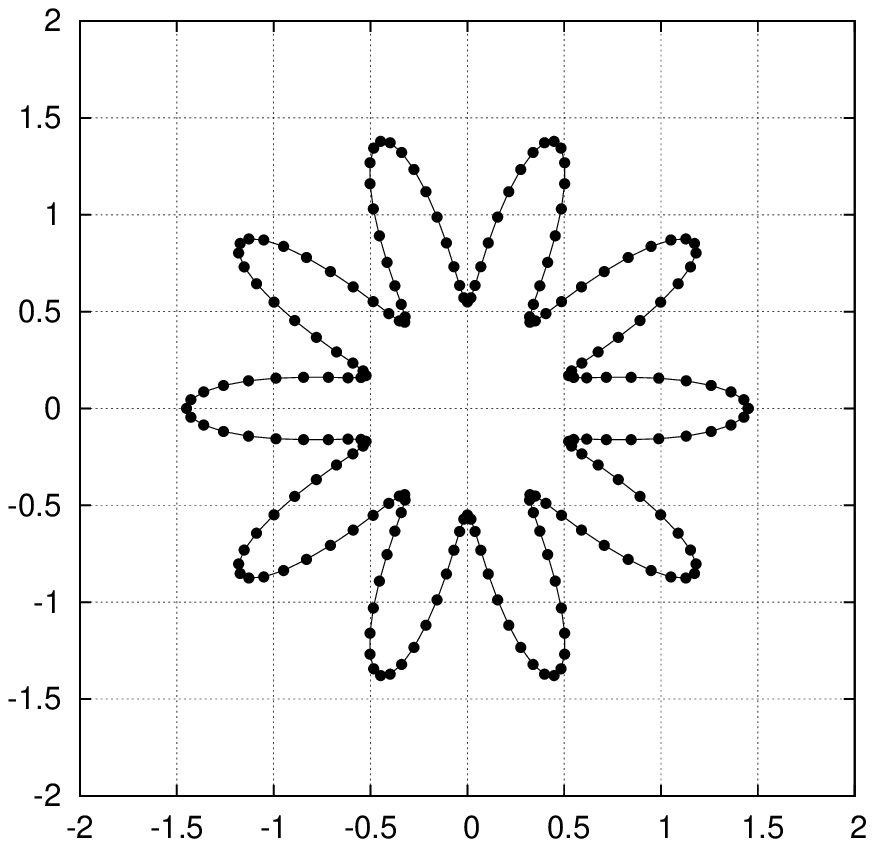}
\includegraphics[width=0.48\textwidth]{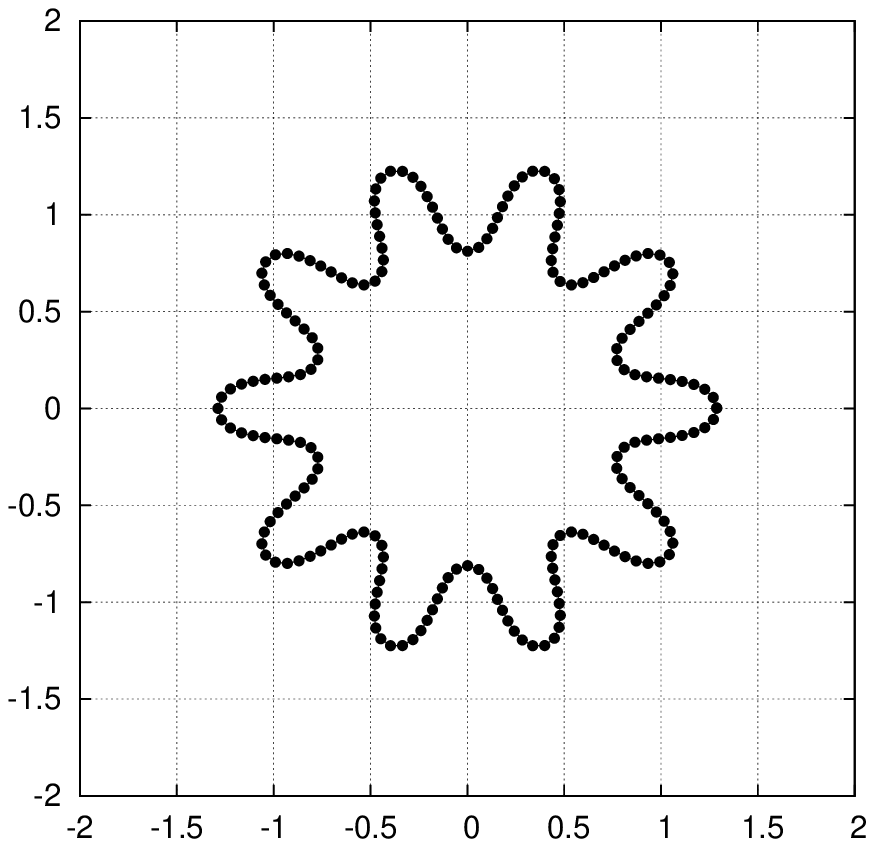}\\
\includegraphics[width=0.48\textwidth]{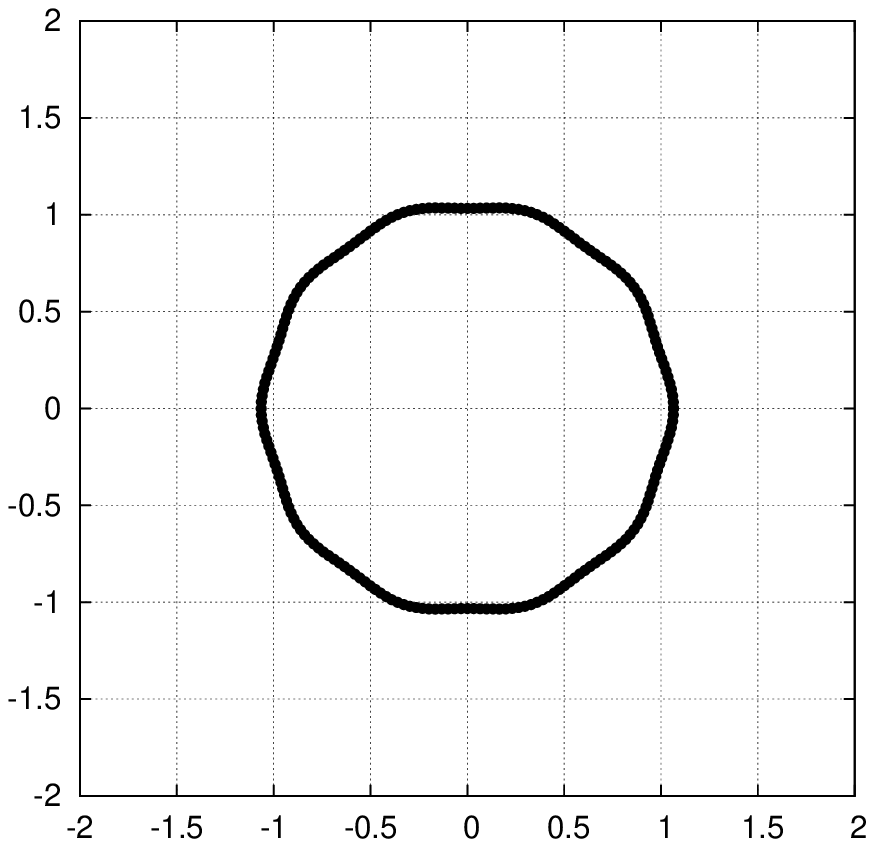}
\includegraphics[width=0.48\textwidth]{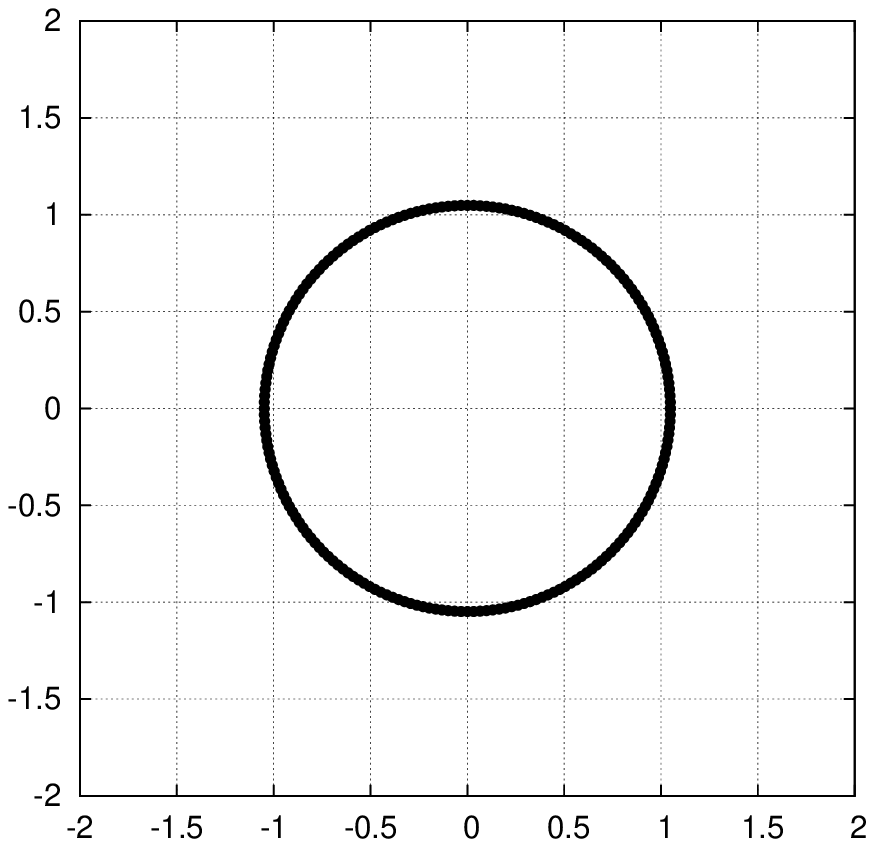}\\
\end{center}
\caption{\rc{Area-preserving} mean curvature flow (\ref{CMCFl}) where the initial 10-folded curve asymptotically attains the circular shape. The curve $\Gamma(t)$ is depicted for $t=0$, $t=0.0125$, $t=0.05$ and $t=0.5$.}
\label{fig:ex3}
\end{figure}

\begin{figure}
\begin{center}
\includegraphics[width=0.48\textwidth]{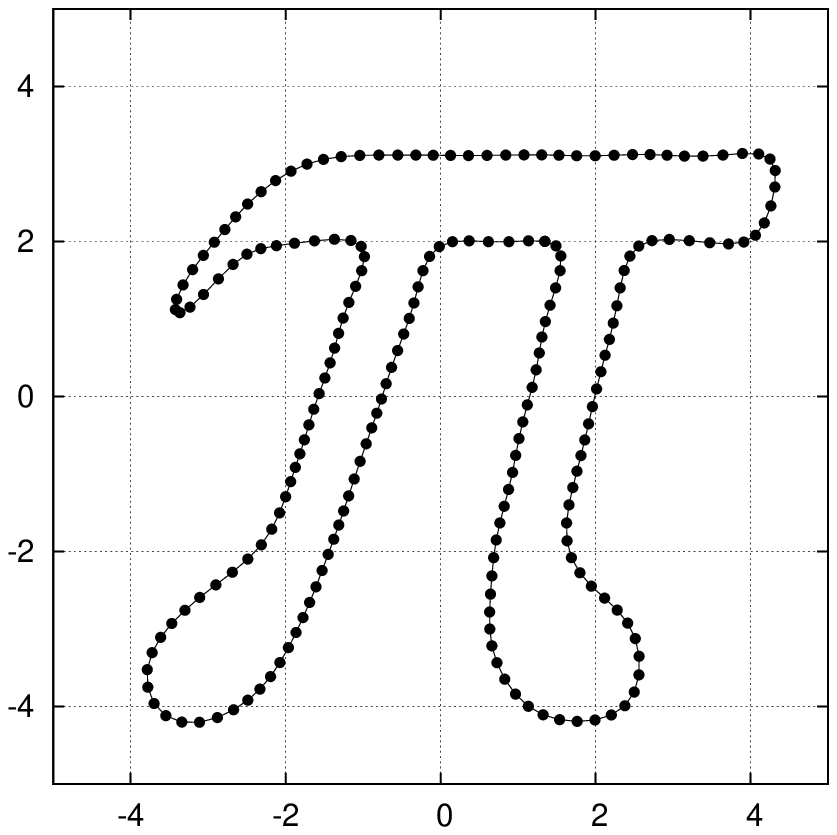}
\includegraphics[width=0.48\textwidth]{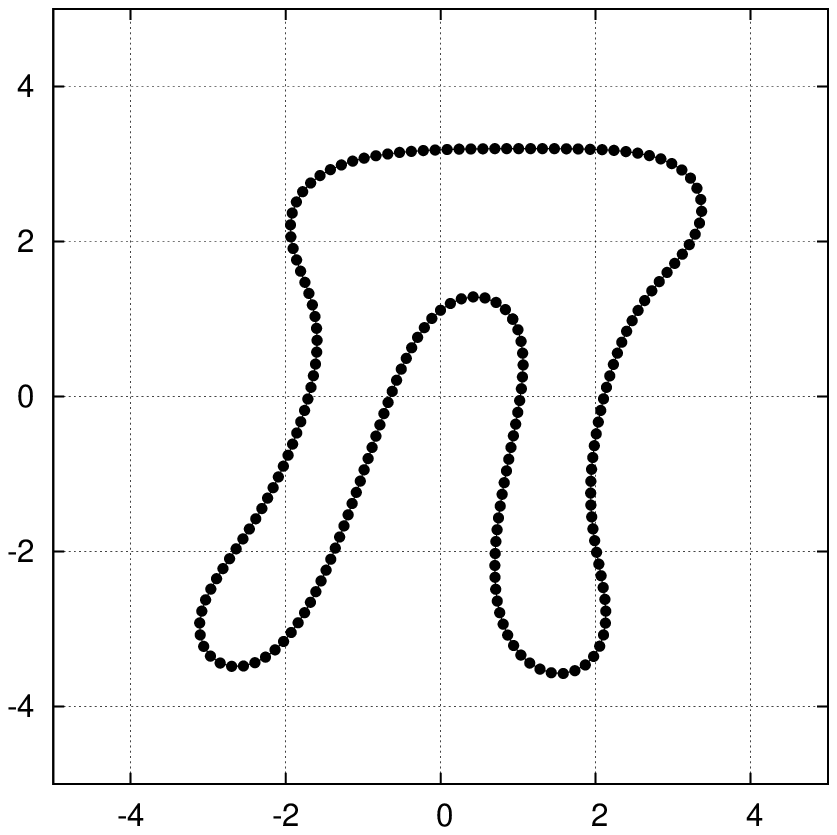}\\
\includegraphics[width=0.48\textwidth]{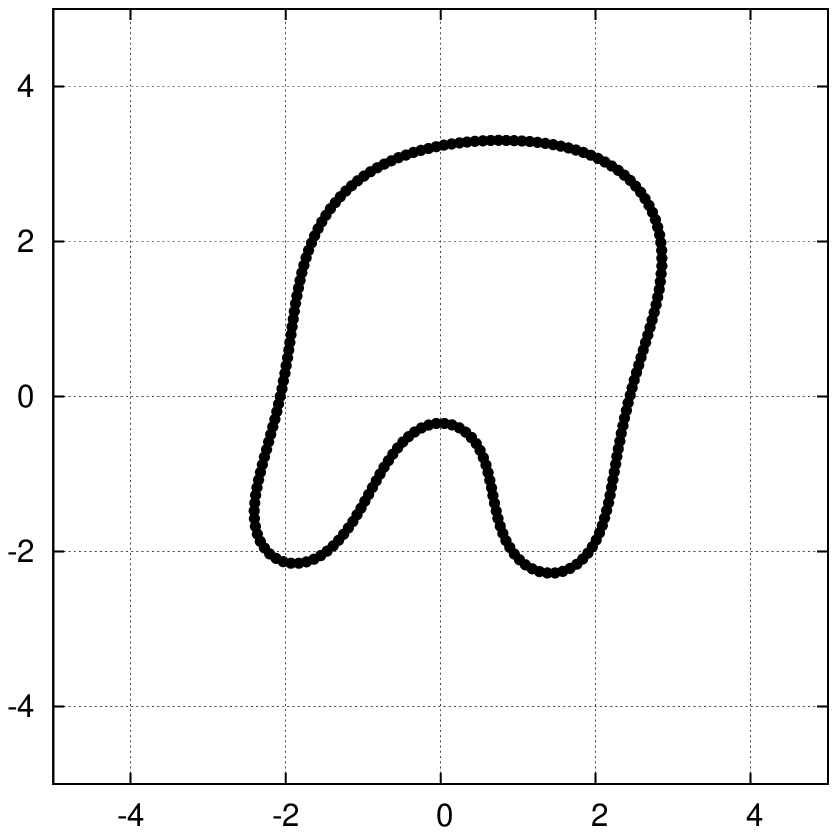}
\includegraphics[width=0.48\textwidth]{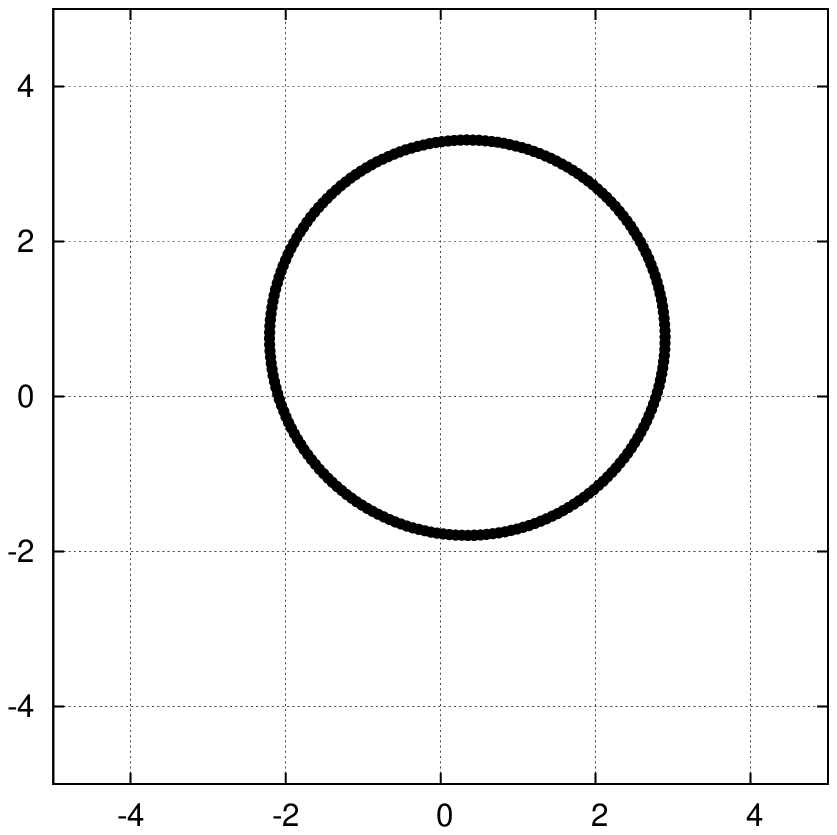}\\
\end{center}
\caption{\rc{Area-preserving} mean curvature flow (\ref{CMCFl}) where the initial $\pi$-shaped curve asymptotically attains the circular shape. The curve $\Gamma(t)$ is depicted for $t=0$, $t=0.05$, $t=0.125$ and $t=1.25$.}
\label{fig:ex4}
\end{figure}


\section{Computational studies}\label{sec4}

We use scheme (\ref{eq:DCiScheme}) to perform a series of computational studies showing the behavior of the solution to (\ref{DirNonl}) as the directly treated constrained mean-curvature flow
$$v_{\Gamma} = - \kappa_{\Gamma}  + \frac{1}{|\Gamma|}\int_\Gamma \kappa_{\Gamma} {\rm d}\rc{s},$$
in comparison to the curve shortening flow:
$$v_{\Gamma} = - \kappa_{\Gamma}.$$
The computations are analyzed using the following measured quantity:
\begin{itemize}
\item Area enclosed by $\Gamma$, 
$A = \int_{{\rm Int} \rc{(\Gamma)}} {\rm d}x$
\rc{should be preserved.}
\end{itemize}

The following examples demonstrate how the solution of (\ref{DirNonl}) evolves in time approaching the circular shape (called the Wulff shape), unlike the usual law (\ref{MCFl}) where the curve shrinks to a point when $F=0$.
\rc{In the examples, the discretization points remain almost uniformly distributed along the evolving curves during the considered evolution time intervals. Any redistribution algorithm was not necessary in this case (compare to \cite{Sinica08,BenJAP10}).}





\noi {\bf Example 1.}
\rc{In Figure \ref{fig:ex1},} the first study shows the behavior of the solution when the initial four-folded curve is given by the formula 
$$ r_0(u) = 1 + 0.4 \cos(8 \pi u), \ u \in [ 0, 1 ].$$
The motion in the time interval $[ 0, 0.5 ]$ is driven by \rc{equation}
(\ref{MCFl}). The curve $\Gamma(t)$ asymptotically approaches the circular shape and shrinks to a point in finite time (compare with \cite{SevcNCMM}, \cite{Gray87}).
The number of finite volumes is $M = 200$.


\noi {\bf Example 2.}
\rc{In Figure \ref{fig:ex2},} the second study shows the behavior of the solution when the initial five-folded curve is given by the formula 
$$ r_0(u) = 1 + 0.65 \cos(10 \pi u), \ u \in [ 0, 1 ].$$
The motion in the time interval $[ 0, 0.5 ]$ is driven by \rc{equation}
(\ref{CMCFl}). The curve $\Gamma(t)$ asymptotically approaches the circular shape whereas the enclosed area is preserved. (see \cite{RuSte92}).
The number of finite volumes is $M = 200$. The initial curve encloses the area of $3.839$ and at $t=0.5$, the curve encloses the area of $3.834$.


\noi {\bf Example 3.}
\rc{In Figure \ref{fig:ex3},} the third study shows the behavior of the solution when the initial ten-folded curve is given by the formula 
$$ r_0(u) = 1 + 0.45 \cos(20 \pi u), \ u \in [ 0, 1 ].$$
The motion in the time interval $[ 0, 0.5 ]$ is driven by \rc{equation}
(\ref{CMCFl}). The curve $\Gamma(t)$ asymptotically approaches the circular shape whereas the enclosed area is \rc{preserved (see \cite{RuSte92}).}
The number of finite volumes is $M = 200$. The initial curve encloses the area of \rc{$3.476$} and at $t=0.5$, the curve encloses the area of \rc{$3.470$}.


\noi {\bf Example 4.}
\rc{In Figure \ref{fig:ex4},} the fourth study shows the behavior of the solution when the initial $\pi$-shaped curve 
{ whose parametric equations can be found in the Wolfram Alpha Database ({\tt http://www.wolframalpha.com}).}
%
%
The motion in the time interval \rc{$[ 0, 1.25 ]$} is driven by problem
(\ref{CMCFl}). The curve $\Gamma(t)$ asymptotically approaches the circular shape whereas the enclosed area is preserved. (see \cite{RuSte92}).
The number of finite volumes is $M = 200$. The initial curve encloses the area of \rc{$20.61$} and at $t=1.25$, the curve encloses the area of $20.53$.


\section{Conclusion}

The paper studies the \rc{area-preserving} mean curvature flow in the terms of qualitative behavior of the solution obtained numerically. The studies confirmed the theoretical indications that the solution approaches the circular shape in long term \rc{(see \cite{Gage86,RuSte92})}. This behavior corresponds to the expected use in modeling the recrystallization phenomena in solid phase.


\section*{Acknowledgement}
The first two authors were partly supported by the \rc{project} No. P108/12/1463 {\it "Two scales discrete-continuum approach to dislocation dynamics"} of the Grant Agency of the Czech Republic \rc{and by the project VEGA 1/0747/12}.




{\small

\bibliographystyle{plain}


\providecommand{\bysame}{\leavevmode\hbox to3em{\hrulefill}\thinspace}
\providecommand{\MR}{\relax\ifhmode\unskip\space\fi MR }
\providecommand{\MRhref}[2]{%
  \href{http://www.ams.org/mathscinet-getitem?mr=#1}{#2}
}
\providecommand{\href}[2]{#2}

\bibliography{phfbase,benes,cmcfl,mach}

}

\vspace{5mm}
\noi{\em Authors' addresses}:
\vspace{5mm}

\noi {\em Miroslav Kol{\' a}{\v r}}, \\
Czech Technical University in Prague, Prague, Czech Republic\\
e-mail: \texttt{kolarmir@fjfi.cvut.cz}.

\noi {\em Michal Bene{\v s}}, \\
Czech Technical University in Prague, Prague, Czech Republic\\
e-mail: \texttt{michal.benes@fjfi.cvut.cz}.

\noi {\em Daniel {\v S}ev{\v c}ovi{\v c}}, \\
Comenius University, Bratislava, Slovakia\\
e-mail: \texttt{sevcovic@fmph.uniba.sk}.

\end{document}